\documentclass[12pt, a4paper]{amsart}
\usepackage[english]{babel}
\usepackage[T1]{fontenc}
\usepackage{amsmath,amsfonts}
\usepackage{amsthm}
\usepackage{amssymb}
\usepackage[utf8]{inputenc}
\usepackage[top=4cm, bottom=4cm, left=3.2cm, right=3.2cm]{geometry}
\title{On the cohomology of elliptic coformal spaces}
\usepackage{systeme}
\usepackage{pifont}
\usepackage[all]{xy}
\usepackage{tikz}
\usepackage{caption}
\usetikzlibrary{fit,calc,positioning,decorations.pathreplacing,matrix}
\newtheorem{thm}{Theorem}[section]
\newtheorem{cor}[thm]{Corollary}

\newtheorem{lem}[thm]{Lemma}

\newtheorem{conj}[thm]{Conjecture}
\newtheorem{notation}[thm]{Notation}

\usepackage{graphicx}
\date{\today}


\begin{document}
\author{Y. Rami}{\let\thefootnote\relax\footnote{{\it Address}: My Ismail University of Meknes, Department of Mathematics  B. P. 11 201 Zitoune,  Meknès, Morocco}\let\thefootnote\relax\footnote{{\it Email}: y.rami@umi.ac.ma}}

\renewcommand{\abstractname}{Abstract}

\selectlanguage{english}

\begin{abstract} {In this article, we pursue the study begun in \cite{Lup02}  on the cohomology of rationally elliptic coformal spaces. Consequently, we complete, for such spaces, the proof of Lupton's conjecture and deduce Hilali's.}
\end{abstract}

\maketitle

\section{Introduction}
Throughout this paper,  vector spaces and  algebras are over the rationals  $\mathbb{Q}$ and, unless otherwise stated, $X$ will be a simply connected CW-complex of finite type. It is said {\it elliptic} if  $\pi_*(X)\otimes \mathbb{Q}$ and  $H_*(X, \mathbb{Q})$ are both finite dimensional and {\it coformal} if the projection  $C_*(\Omega X;\mathbb{Q})\rightarrow (H_*(\Omega X;\mathbb{Q}),0)$ is a quasi-isomorphism of differential graded algebras. There are two practical equivalent characterizations to coformality. The first states that, $(L_X = \pi_*(X)\otimes \mathbb{Q},0)$, the rational homotopy Lie algebra equipped with the zero differential, is a {\it Lie model} of $X$ ($\S 2$) and the second stipulates that its {\it minimal Sullivan model} $(\Lambda V,d)$ (or {\it model} for short) has a quadratic differential ($\S 2$).
 
 One of the first invariants introduced in the study of homotopy properties of topological spaces is the {\it Lusternik-Schnirelmann category} ({\it LS-category} for short). For an arbitrary topological space $X$, this is  denoted $cat(X)$ and defined  as  the smallest number $n$ such that $X$ is covered by $n+1$  contractible open sets. Being difficult to calculate, it was approximated by different invariants of an algebraic nature. One of its  closer lower bounds is the {\it rational LS-category} $cat_0(X)$ defined as the LS-category $cat(X_0)$ of the rationalization  $X_0$  of $X$. Using the minimal models, Y. Félix and S. Halperin gave an algebraic caracterization of $cat_0(X)$ as being the least integer $m$ such that the projection $(\Lambda V,d) \stackrel{pr_m}{\rightarrow} (\Lambda V/\Lambda^{\geq m+1}V,\bar{d})$ admits a retraction as a morphism of  differential graded algebras \cite{FH82}.  This is furthermore lowered by the {\it rational Toomer invariant} denoted $e_0(X)$ and defined in terms of the {\it Milnor-Moore spectal sequence}: 
 \begin{equation}\label{TopM-Mss}
  Ext^{p,q}_{H_*(\Omega X, \mathbb{Q})}(\mathbb{Q , \mathbb{Q}}) \Rightarrow H^{p+q}(X , \mathbb{Q}).
  \end{equation}
 as being the smallest integer $p$ such that the  $E_{\infty}^{*}$ is zero starting from $p+1$.  For elliptic coformal spaces, thanks to Poincar\'e duality property, $e_0(X)$ and $cat_0(X)$ coincide and are both equal to   $\dim (\pi_{*}(X)\otimes \mathbb{Q})_{odd}$ \cite[Proposition 10.6]{FH82}.
 
The emphasis on the study of coformal spaces follows the equivalence, established in \cite[Proposition 9.1]{FH82} between the above spectral sequence and its {\it algebraic version}:
\begin{equation}\label{AlgM-Mss}
H^{p,q}(\Lambda V,d_2) \Rightarrow H^{p+q}(\Lambda V,d).
\end{equation}
where $d_2 : V \rightarrow \Lambda^2 V$ is  the quadratic part  of $d$. In particular, any property of the cohomology of $(\Lambda V,d_2)$ (which is coformal) will have an impact on that of $(\Lambda V,d)$ (see for instance \cite[Theorem 1.1]{Ram19}).

Recall that, for any coformal space $X$ with model $(\Lambda V,d)$, the lengths of representative cocycles come with a lower graduation on $H(\Lambda V,d)$ and hence on  $H^*(X, \mathbb{Q})$  so  that, 
$$H^l(X, \mathbb{Q})=\oplus_{k\geq 0}H^l_k(X, \mathbb{Q}),\; \hbox{for each} \; l\geq 0.$$
   
 Recall also that, for any topological space $Y$, the homology classes in $H_*(Y ; \mathbb{Q})$ lying in the image of the
 Hurewicz homomorphism $hur_Y: \pi_*(Y)\otimes \mathbb{Q}\rightarrow H_*(Y, \mathbb{Q})$ are called {\it spherical classes}. Dually, a cohomology class $\alpha \in H^k(Y ; \mathbb{Q})$ is {\it spherical}
 if there exits $f: \mathbb{S}^n \rightarrow Y$ such that $f^*(\alpha)\not =0$ (see below for more details).
  
 Our main result   is stated as follows:
\begin{thm} \label{Thm1}
Let $X$ be a coformal elliptic space whose graded algebra $H^*(X, \mathbb{Q})$ is generated by at least two generators, and all its spherical classes have even degrees. Then,
\begin{equation}\label{dim}
\dim H^*_k(X, \mathbb{Q})\;   \;  \left\{ \begin{array}{ll}
= 1\; if \; k = 0\; or \;  e_0(X)\\
\geq 2\; if \; 1\leq k \leq e_0(X)-1\\
= 0\; if\;  k\geq e_0(X)+1.
\end{array} \right.
\end{equation}
\end{thm} 


In \cite[Theorem 2.5]{Lup02} G. Lupton established, subject to the presence in $H^*(X, \mathbb{Q})$ of at least one spherical class of odd degree, the property  (\ref{dim}) for rational spaces having a model $(\Lambda V, d)$ with a homogeneous differential of length $l\geq 2$. He then conjectured that this is generally true for such spaces without additional assumptions (cf. \S 2 below). Consequently, we have:
\begin{thm}\label{Lup-Top-conj}
The Lupton conjecture is satisfied for any elliptic coformal space.
\end{thm}
In \cite{Hil90}, M.R. Hilali posed a rather old conjecture which states that {\it for any elliptic space $Y$ we have $\dim H^*(Y,\mathbb{Q}) \geq \dim \pi_*(Y)\otimes \mathbb{Q}$}.  
In \cite{Ram19},
 we 
showed (cf. Theorem 1.2) that Hilali's conjecture  is true, in particular, for any coformal space whose rational homotopy is concentreted in odd degrees. 
As a second consequence of the Theorem \ref{Thm1}, we enhance \cite[Theorem 1.2]{Ram19} by supressing the additional hypothesis:
\begin{thm}\label{Hil-Top-conj}
The Hilali conjecture is satisfied  for any elliptic  coformal space. 
\end{thm}
The proof of this theorem is given in terms  of  its algebraic-version, namely,  Corollary 3.5 to follow.
 
As examples of spaces  satisfying the hypothesis  of Theorem \ref{Thm1}, we cite the class of $F_0$-spaces which are characterized by their  rational cohomology  of the form
 $\mathbb{Q}[x_1, \ldots x_n]/(f_1, \ldots , f_n)$ where each $f_i$ is a nonzero homogeneous polynomial of degree two.  Among these, we have  $\mathbb{C} P^2 \# \mathbb{C} P^2$ \cite[Example 7.3]{FH82}.
 Other examples far from the latter class are the homogeneous spaces $Sp(6)/SU(5)$  and $Sp(6)/SU(3)\times SU(3)$  \cite{Mur94}. Likewise,  referring to  \cite[Exemple 2]{Hua}, the total space $E$ of the fibration $\mathbb{S}^{2k+1}\hookrightarrow E \rightarrow \vee_{i=1}^{i=n}\mathbb{S}^{2n_i}$ ($k\geq 1$) with rationally nontrivial inclusion $\mathbb{S}^{2k+1}\hookrightarrow E$  is an elliptic coformal space which moreover satisfies hypothesis of  Theorem \ref{Thm1}. 
\section{Main tools}
Recall from the introduction that $X$ stands for a finite-type simply connected CW-complex and the ground field is $\mathbb{Q}$.
\subsection{Sullivan minimal models}
Let $V= \oplus_{i\geq 0}V^i$ be a graded vector space  where $V^i$ denotes the subspace of elements $v\in V$ of homogeneous degree $|v|=:i$. The commutative graded free algebra, denoted $\Lambda V$, on $V$ is the quotient of the free graded algebra $TV$ with the graded ideal generated by elements of the form $v\otimes v' - (-1)^{|v||v'|} v'\otimes v$. That is:
$$\Lambda V = Exterior(V^{odd}) \otimes Symmetric(V^{even})$$ where $V^{even} = \oplus _{i\geq 0}V^{2i}$ and $V^{odd} = \oplus _{i\geq 0}V^{2i+1}.$ 

Assume that $V$  has a well-ordered   basis $\{x_{\alpha}\}$ satisfying $dx_{\alpha }\in \Lambda V_{<\alpha}$  where $V_{<\alpha} = \{v_{\beta} \mid \beta < \alpha \}$ and let $d : V \rightarrow \Lambda V$  a linear map of degree $+1$. This is extended to  a derivation $d : \Lambda V \rightarrow \Lambda V$ by putting:
$$d(x_{\alpha}x_{\beta}) = d(x_{\alpha})x_{\beta} + (-1)^{|x_{\alpha}|}x_{\alpha}d(x_{\beta}).$$
If morover $d^2=0$, $(\Lambda V, d)$ becomes a free commutative differential graded algebra. It is said a {\it Sullivan algebra}. It is said a {\it minimal Sullivan algebra} if in addition $deg(x_{\alpha })< deg(x_{\beta })$ implies $\alpha
 <\beta $. When  $V^0=\mathbb{Q}$ and  $V^1=0$, minimality  is equivalent to decomposability of $d$ in the sens that $ d(V)\subseteq \oplus _{i\geq 2} \Lambda ^iV  =: \Lambda ^{\geq 2}V$. 
 
 Recall from Sullivan theory  that there is a unique (up to isomorphism) minimal Sullivan algebra  $(\Lambda V,d)$ and a quasi-isomorphism, i.e. a morphism inducing an isomorphism in cohomology, $m_X : (\Lambda V,d) \stackrel{\simeq }{\rightarrow} A_{PL}(X)$ with source the algebra of piecewise-linear de Rham forms on $X$ \cite{Sul78}.   By hypothesis on $X$, we have $V^0=\mathbb{Q}$ and  $V^1=0$ and therefore $m_X$ or simply $(\Lambda V,d)$ is called the {\it minimal Sullivan model} (or model for short) of $X$.  Moreover, $X$ and $(\Lambda V,d)$ are linked as follows:
 \begin{equation}
V^i\cong Hom_{\mathbb{Z}} (\pi_i(X), \mathbb{Q}),\; (i\geq 2)\; \hbox{and}\; H^*(\Lambda V,d)\cong H^*(X,\mathbb{Q}).
 \end{equation} 
 
The aim of this article is to continue the in-depth study on the cohomology of coformal elliptic spaces begun by G. Lupton in \cite{Lup02}.
 Recall from the introduction that, in terms of its model $(\Lambda V,d)$, $X$ is {\it elliptic} if and only if $\dim V < \infty$ and $\dim H(\Lambda V,d) < \infty$. It is {\it coformal} if the differential $d$ satisfies $d: V \rightarrow \Lambda ^{2} V$ or, equivalently, it is of homogeneous-length $2$. There is then a cochain complex:
 \begin{equation}
 \ldots  \rightarrow \Lambda^{k-1}V \stackrel{d}{\rightarrow} \Lambda^{k}V \stackrel{d}{\rightarrow} \Lambda^{k+1}V \rightarrow \ldots 
 \end{equation}
  which induces on cohomology a lower  graduation given by lengths of cocycle representatives: 
  \begin{equation}\label{decomposition}
 H^*(\Lambda V,d)=\oplus_{k\geq 0}H^*_k(\Lambda V,d).
  \end{equation}
 \subsection{Sphyrical cohomology classes}  
Recall also from the introduction that $\alpha \in \tilde{H}^r(X, \mathbb{Q})$  is {\it spherical} if there exists $[f]\in \pi_r(X)\otimes \mathbb{Q}$ such that $f^*(\alpha)\not = 0$ where $f^* : \tilde{H}^r(X, \mathbb{Q})\rightarrow \tilde{H}^r(\mathbb{S}^r, \mathbb{Q})$ is the induced morphism by $f$. If  $[a]\in  \tilde{H}_r(\mathbb{S}^r, \mathbb{Q})$ is the generating class and $f_* : \tilde{H}_r(\mathbb{S}^r, \mathbb{Q})\rightarrow \tilde{H}^r(X, \mathbb{Q})$, we know that  $hur_X([f]) = f_*([a])$. It results that 
$<f^*(\alpha) , [a]> = <\alpha , hur_X([f])> = <hur_X^*(\alpha) , [f]>$;
$hur_X^*$ being the dual of Hurewisz morphism $hur_X: \pi_{*}(X)\otimes \mathbb{Q}\rightarrow \tilde{H}_*(X,\mathbb{Q})$.
Therefore, $\alpha$ is  spherical  if and only if $hur_X^*(\alpha)\not = 0$. 

Now, if $(\Lambda V ,d)\stackrel{m_X}{\rightarrow} A_{PL}(X)$ is a model of $X$, referring to $\S 13(c)$ in \cite{FHT01} we have an identification between $hur_X^*$ and the projection $\zeta : H^+(\Lambda V ,d) \rightarrow V\cap \ker(d)$ given by 
$$<\zeta ([z]) , [f]> = <H(m_X)([z]) , hur_X([f])>, \; \hbox{for any}\; [z]\in H^+(\Lambda V ,d).$$ Thus for any cohomology class $\alpha \in H^*(X , \mathbb{Q})$ and $[z] = H(m_X)^{-1}(\alpha)$  we deduce that $\alpha$ is  spherical  if and only if $\zeta([z])\not = 0$. As a conclusion, { spherical cohomology classes are determined by some elements} in $V\cap \ker(d)$.
\subsection{Rational Toomer invariant}
 Given a rational elliptic space $X$. Referring to \cite{FHT91} we know that it satisfies {\it Poincaré duality property}  in the sens that  for some integer  $N$, $H^{>N}(X, \mathbb{Q})=0$, $H^{N}(X, \mathbb{Q})\cong \mathbb{Q}$ and, if $\omega \in H^{N}(X, \mathbb{Q})$ is a generating class then $\cap \omega : H^{i}(X, \mathbb{Q}) \rightarrow H_{N-i}(X, \mathbb{Q})$ is an isomorphism for all $0\leq i \leq N$. $\omega$ is called the  {\it fundamental class} of $X$ and $N$ its {\it formal dimension}. An essential tool attached to such spaces is the rational Toomer invariant defined as follows (see for instance \cite{FH82}):
$$e_0= sup\{ k\;    \hbox{such that} \; E_{\infty}^{k,*}\not =0 \}$$
where $E_{\infty}^{k,*}$ stands for the $\infty$ term of the Milnor-Moore spectral sequence (\ref{TopM-Mss}). Once more, the equivalence between (\ref{TopM-Mss}) and (\ref{AlgM-Mss}) gives us the following convenient formula to determine or at least to approximate this invariant quite easily:
$$e_{0} = sup\{k\; | \; \omega \; \hbox{can be represented by a cocycle in} \; \Lambda^{\geq k}V\}.$$
Now, for any $x\in H^*(\Lambda V,d)$ one can define  its {\it Toomer invariant}:
 $$e_0(x) =sup\{k\; | \; x \; \hbox{can be represented by a cocycle in} \; \Lambda^{\geq k}V\}$$ so that $e_0 = e_0(\omega)$ \cite{FH82}.
\section{Algebraic statements of our results}
In his famous article \cite{Lup02}, G. Lupton showed that for any  elliptic Sullivan model $(\Lambda V,d)$ with homogeneous differential of constant length $l\geq 2$, we have 
$$H_k(\Lambda V,d)\not =0,\; \forall \; 0\leq k \leq e_0.$$ Then, using a more profound analysis of the cohomology of such models, he established the following 
\begin{thm}\label{LuTh}\cite[Theorem 2.5]{Lup02}
 Suppose  $(\Lambda V,d)$ is an elliptic Sullivan algebra with a homogeneous differential of length $l\geq 2$ and $\ker(d: V^{odd}\rightarrow \Lambda V)$ is non-zero. Then $\dim  H^*_k(\Lambda V,d)\geq 2$ for each $k= 1, \ldots ,e-1$, where $e=e_0(\Lambda V,d) = \dim V^{odd}+(l-2)\dim V^{even}$.
 \end{thm} 
This causes him to pose the following  conjecture which we call henceforth   {\it Lupton's conjecture}
\begin{conj}\label{Conjecture}
Let $(\Lambda V,d)$ be an elliptic Sullivan algebra with a homogeneous differential of length $l\geq 2$. Either $\dim H^*_k(\Lambda V,d)\geq 2$ for $k=1, \ldots , e-1$  or $H^*(\Lambda V,d)$ is a truncated polynomial algebra
on a single generator.
\end{conj}
Notice that under the hypothesis of the above theorem, $H^*(\Lambda V,d)$ can not have the structure of a truncated polynomial algebra on a single generator. Moreover, according to the above section, we see that the degree of any spherical cohomology classes  $[z]\in H^+(\Lambda V,d)$ is exactly that of $\zeta([z])\in V\cap \ker(d)$. Therefore, the above theorem gives the proof of conjecture (\ref{Conjecture}) in the presence of at least one spherical cohomology class of odd degrees. 

In this article, we limit ourselves to coformal elliptic Sullivan models, that is to say when the length of the differential is $l=2$. Under the complementary hypothesis, namely, when all the of spherical cohomology classes are of even degree, we establish the following algebraic translation of Theorem \ref{Thm1}.
\begin{thm}\label{alg-ver}
 Let $(\Lambda V,d)$  be an elliptic coformal Sullivan algebra. If
 \begin{enumerate}
 \item[(a)] $\ker(d: V\rightarrow \Lambda V)\subseteq V^{enen}$,
 \item[(b)] the commutative graded algebra $H^*(\Lambda V,d)$ is generated by at least two generators,
 \end{enumerate}
 then, $\dim H^*_k(\Lambda V,d)\;  is \;  \left\{ \begin{array}{ll}
 \geq 2,\; if \ 1\leq k \leq \dim V^{odd}-1\\
 = 0,\; if\;  k\geq \dim V^{odd}+1.
 \end{array} \right.$
  \end{thm} 
Remark that the cohomology $H^*(\Lambda V,d)$ of a coformal Sullivan model $(\Lambda V,d)$ has the structure of truncated polynomial algebra if and only if $(\Lambda V,d)=(\Lambda (x_1, x_2), d)$ where $x_1$ has an even degree and $d(x_2) = x_1^2$, that is to say, if and only if  $(\Lambda V,d)$ is the model of an even sphere (cf. Lemma 4.3 below).

 According to the discussion above, we obtain the complete algebraic-version of Theorem \ref{Lup-Top-conj} as follows:
\begin{cor}
Lupton's conjecture is satisfied for any coformal elliptic Sullivan algebra $(\Lambda V , d)$. 
\end{cor}
Next, using Theorem \ref{alg-ver}, 
combined with \cite[Corollary 2.6]{Lup02}, 
and the fact that for any elliptic space, $\dim V^{odd}\geq \dim V^{even}$, we  get  an affirmative answer to  the following algebraic version of Corollary \ref{Hil-Top-conj}:
\begin{cor}
Hilali's conjecture is satisfied for any elliptic coformal Sullivan algebra ($\Lambda V,d_2)$.
\end{cor}
Let us  note in passing that this enhances the result obtained by Ben El Krafi et al. in \cite{EHM}. 
\section{Proof of Theorem \ref{alg-ver}}
This section is devoted to the proof of our main result. We began by proving some lemmas which permit us to emphasis on the essential part of the proof.
\subsection{The bigraded Gysing exact sequence and preparatory results}, is unique (up to isomorphism). It is  closely related to $X$ by the isomorphisms $V^i\cong Hom_{\mathbb{Z}} (\pi_i(X), \mathbb{Q}),\; (i\geq 2)\; \hbox{and}\; H^*(\Lambda V,d)\cong H^*(X,\mathbb{Q})$.
In this section, we present specific tools that will be used to prove our main result. Recall that  we are concerned with coformal elliptic spaces each of which is endowed with a model $(\Lambda V,d)$  satisfying:
$$\dim V <\infty,\; \dim H^*(\Lambda V, d) <\infty\; \hbox{and}\; d(V)\subseteq \Lambda ^2V.$$
For the remainder, we  denote
\begin{equation}
(\Lambda V,d)=: \Lambda (x_1, x_2, \ldots ,x_n ,d) \; \hbox{with} \; |x_1|\leq |x_2|\leq \ldots \leq |x_n|.
\end{equation} 
\begin{notation}
With the notations above, there is a short exact sequence of differential graded vector spaces
\begin{equation}
0\rightarrow (x_1\wedge V, d)\stackrel{\subset}{\rightarrow}(\Lambda V,d)\stackrel{p}{\rightarrow} (\Lambda W, \bar{d})\rightarrow 0
\end{equation}
where  $x_1\wedge V$ is the ideal of $\Lambda V$ generated by $x_1$, $W= \mathbb{Q}x_2\oplus \mathbb{Q}x_3 \oplus \cdots \oplus \mathbb{Q}x_n$ and $\bar{d}$ the differential deduced from the isomorphism of graded algebras $\Lambda W \cong \Lambda V/x_1\wedge  V$.
\end{notation}
\begin{lem}
The above exact sequence induces the following one
\begin{equation}\label{Gysin} 
\ldots \stackrel{\delta ^*}{\rightarrow} H_{k-1}^{i-2r}(\Lambda V,d)\stackrel{j^*}{\rightarrow}H_{k}^{i}(\Lambda V,d)\stackrel{p^*}{\rightarrow}H_{k}^{i}(\Lambda W,\bar{d})\stackrel{\delta ^*}{\rightarrow}H_{k}^{i-2r+1}(\Lambda V,d)\rightarrow \ldots
\end{equation}
called the {\it Gysin sequence}.
\end{lem}
\begin{proof}
Denote by $j: \Lambda V\rightarrow \Lambda V$  the map of degree $|j|= 2r$ defined by $j(\chi)=x_1\chi,\; \hbox{for any }\; \chi \in \Lambda V$ and consider the following short  exact sequence of differential graded vector spaces:
$$0\rightarrow \Lambda V\stackrel{j}{\rightarrow} \Lambda V \stackrel{p}{\rightarrow}\Lambda W\rightarrow 0.$$ The induced long exact sequence in cohomology is neither than (\ref{Gysin}).  Its connecting morphism is defined for any class  $[\chi ]\in H_{k}^{i}(\Lambda W,\bar{d})$ by $\delta ^*([\chi ]) = [\chi ']$ with $\chi '\in \Lambda V$ is such that $d\chi = x_1\chi ' = j(\chi ')$ or equivalently  $\bar{d}\chi  = 0$ 
 \cite[Proof of Theorem 2.2 (case II)]{Lup02}.
 \end{proof}
We assume once and for all that $(\Lambda V, d)$  is an elliptic Sullivan model satisfying hypothesis
of Theorem \ref{alg-ver} with:
\begin{equation}\label{cond(V)}
V = \mathbb{Q}x_1 \oplus W; \; W =\mathbb{Q}x_2  \oplus \mathbb{Q}x_3  \oplus \ldots \oplus \mathbb{Q}x_n.
\end{equation}
It results, for degree reason and by hypothesis (a) of that theorem that$dx_1=0$ and $|x_1| =2r$ (some $r\geq 1$) is even.

To make the proof of Theorem \ref{alg-ver} clearer, we first establish the following  preparatory lemmas:
\begin{lem}\label{dimV}
Let $(\Lambda V,d)$ be an elliptic coformal Sullivan algebra satisfying the hypothesis of Theorem \ref{alg-ver}.
 Then the graded vector space $V$ satisfies:
 \begin{enumerate}
\item $\dim V\geq 4$.
\item $\dim V^{odd}\geq 3$ or else $(\Lambda V,d)$ is quasi-isomorphic to $(\Lambda (x_1, x_2, x_3, x_4),d)$ with $d(x_1)=d(x_2)=0$, $d(x_3)=x_1^2$ and  $d(x_4)=x_2^2$  (i.e. a minimal Sullivan algebra of $\mathbb{S}^{2q}\times \mathbb{S}^{2q'}$ some $q, q'\geq 1$).
 \end{enumerate}
\end{lem}
\begin{proof}
First, observe that hypothesis (a) of  Theorem \ref{alg-ver} implies $\dim V \geq 1$. 
\begin{enumerate}
\item 

If $\dim V = 1$ the degree of $x_1$  should be odd which contradicts hypothesis (a). 

If $\dim V =2$, the ellipticity of $(\Lambda V,d)$ implies  $d(x_2)=x_1^2$ and $H^*(\Lambda V, d)\cong \mathbb{Q}[x_1]/(x_1^2)$ is generated by only one generator and  (b) is then not satisfied. 

If $\dim V =3$, since, by ellipticity   $\dim V^{odd}\geq \dim V^{even}$, we should have  $\mid x_2\mid$ and $\mid x_3\mid$ are both odd, so one of them is necessarily a cocycle,  that is  (a) is again not satisfied.

It results that the cases where $\dim V = 1,\; 2,\; 3$ are ruled out and, consequently,   $\dim V\geq 4$.
\item Since $\dim V^{odd}\geq \dim V^{even}$ and $\dim V\geq 4$,  $\dim V^{odd}= 1$ is excluded. Now, if $\dim V^{odd}=2$ necessarily $\dim V^{even}=2$ hence, the only option is the one described in the statement.
\end{enumerate}
\end{proof}
\begin{lem}\label{LC(4)}
Every elliptic coformal Sullivan algebra $(\Lambda V,d)$ with $\dim V \leq 4$ verifies Lupton's conjecture.
\end{lem}
\begin{proof} 
 If $\dim V = 1\; \hbox{or}\; 3$, the discussion made in the proof of the above lemma shows that one of the generating elements of $V$ is a cocycle of odd degree, that is we are in the condition of Theorem \ref{LuTh}, so Lupton's conjecture is verified. Now, if $\dim V =2$, then (cf. the same discussion) $H^*(\Lambda V,d)\cong \frac{\mathbb{Q}[x_1]}{(x_1^2)}$ is a truncated algebra and again Lupton's conjecture is verified. Next, suppose that $\dim V =4$. We have two cases:\\
 - If $\ker(d: V^{odd}\rightarrow \Lambda V)$ is non-zero, we use Theorem \ref{LuTh} to conclude.\\
 - If $\ker(d: V\rightarrow \Lambda V)\subseteq V^{enen}$, then by coformality,  $H^*(\Lambda V,d)$  is a truncated polynomial algebra if and only if $(\Lambda V,d)$ is quasi-isomorphic to $(\Lambda (x_1,x_2), d))$ with $|x_1|$ even and $dx_2=x_1^2$, that is if and only if $\dim V =2$. Thus, the hypothesis $\dim V =4$ puts us in the conditions of the above lemma. Now, if $\dim V^{odd} =3$, by minimality, $(\Lambda V,d)$ should be quasi-isomorphic to $(\Lambda (x_1, x_2, x_3, x_4),d)$ with $|x_1|$  even, $dx_1=0$, $|x_2|\leq |x_3|\leq |x_4|$ are all odd,  $d(x_2)= \alpha_1x_1^2\not =0$,  $d(x_3)= \alpha_2x_1^2\not =0$ and  $d(x_4)= \alpha_3x_1^2 + \alpha_4 x_2x_3\not =0$. But, $\frac{1}{\alpha_2}x_3 - \frac{1}{\alpha_1}x_2$ is clearly a cocycle of odd degree which contradicts the assumption $\ker(d: V\rightarrow \Lambda V)\subseteq V^{even}$. We then conclude by using the assertion (2) of the above lemma.
 \end{proof}
\subsection{Proof of Theorem \ref{alg-ver}}
To be more clear, we divided the proof into three steps: $\dim H^*_1(\Lambda V,d)\geq 2$, $\dim H^*_2(\Lambda V,d)\geq 2$  and $\dim H^*_k(\Lambda V,d)\geq 2$ for $k\geq 2$.  Now, under our hypothesis, Lemma \ref{dimV} and Lemma \ref{LC(4)} require respectively to assume, from now on, that $n=\dim V\geq 4$  and that Theorem \ref{alg-ver} is satisfied for any $(\Lambda W, d)$ with $3\leq \dim W \leq n-1$. This implies that $H^*(\Lambda W, d)$ is not a truncated polynomial algebra on a single generator.
\begin{notation}

 In all what follows, $\mathbb{Q}.[x _i]$ will denote the one-dimensional vector space generated by the cohomology class $[x_i]$ whereas $\mathbb{Q}[x_i]$ will designate the usual polynomial algebra on one determinate $x_i$.
\end{notation}
By hypothesis $(a)$ of the theorem, $|x_1|=2r\geq 2$. \underline{\it  In the sequel,   we put $|x_2| = m_1$.}
\subsubsection{First step: $\dim H^*_1(\Lambda V, d) \geq 2$}
We  consider the following exact sequence which comes from (\ref{Gysin})  when $k=1$:
\begin{equation}\label{Gysin(1)} 
\ldots \stackrel{\delta ^*}{\rightarrow} H_{0}^{i-2r}(\Lambda V,d)\stackrel{j^*}{\rightarrow}H_{1}^{i}(\Lambda V,d)\stackrel{p^*}{\rightarrow}H_{1}^{i}(\Lambda W,\bar{d})\stackrel{\delta ^*}{\rightarrow}H_{1}^{i-2r+1}(\Lambda V,d)\rightarrow \ldots
\end{equation}
 Notice that since $X$ is assumed simply connected then $V^1 =0$. We first put $i=2r$ so that (\ref{Gysin(1)}) induces the following short exact sequence
\begin{equation}\label{Gysin(1) (2r)}
0\rightarrow H_0^{0}(\Lambda V,d)=\mathbb{Q}\rightarrow H_1^{2r}(\Lambda V,d)\stackrel{p^*}{\rightarrow} H_1^{2r}(\Lambda W,\bar{d})\rightarrow 0,
\end{equation}
which means that $\ker(p ^*) =Im(j^*)=\mathbb{Q}.[x_1]$. Hence,
\begin{equation}\label{(a)}
H_1^{2r}(\Lambda V,d) \cong \mathbb{Q}.[x_1] \oplus H_1^{2r}(\Lambda W,\bar{d}).
\end{equation} We need then to consider two cases:
\subsubsection*{Assume $m_1=2r$}
This is equivalent to say that $|x_2| = |x_1|$ so, using the induction hypothesis for $(\Lambda W,\bar{d}$, we have $H_1^{2r}(\Lambda W,\bar{d})\not =0$ which by (\ref{(a)}) give us $\dim H^*_1(\Lambda V,d)\geq 2$.
\subsubsection*{Assume $m_1>2r$} That is we assume $H_1^{2r}(\Lambda W,\bar{d})=0$ so that $\dim H_1^{2r}(\Lambda V,d) =1$. We then continue our cheeking by reconsidering  (\ref{Gysin(1)}) with  $i=m_1$.  This leads us  to  the exact sequence:
\begin{equation}\label{Gysin(1) (m_1)}
0\rightarrow H^{m_1}_1(\Lambda V,d)\stackrel{p^*}{\rightarrow} H^{m_1}_1(\Lambda W,\bar{d})\stackrel{\delta ^*}{\rightarrow}H^{m_1-2r+1}_1(\Lambda V,d)\stackrel{j^*}{\rightarrow}H_2^{m_1+1}(\Lambda V,d)\rightarrow \ldots
\end{equation}
Using  minimality of the Sullivan algebra $(\Lambda W,\bar{d})$ we deduce that $[x_2]\in H^{m_1}_1(\Lambda W,\bar{d})$ is non-zero. If $\dim H^{m_1}_1(\Lambda W,\bar{d})\geq 2$, then  $\dim H^{m_1-2r +1}_1(\Lambda V,{d}) + \dim H^{m_1}_1(\Lambda V,d)\geq 2$ and we are done. Orherwise, $\dim H^{m_1}_1(\Lambda W,\bar{d})=1$, so, by (\ref{Gysin(1) (m_1)}), either $H^{m_1}_1(\Lambda V,d)\not =0$ which permit to conclude or else $H^{m_1-2r+1}_1(\Lambda V,d)\not =0$ and we have to consider the case where $m_1-2r+1=2r$ that is where $m_1=4r-1$.  Now, if $m_1=4r-1$, which is odd,  using Theorem \ref{LuTh}  for $(\Lambda W,\bar{d})$ we get  $\dim H^{m}_1(\Lambda W,\bar{d}) \geq 1$ for some integer $m > m_1$. We then use once again (\ref{Gysin(1) (m_1)})  with $m$ replacing $m_1$ which leads to  
 $H^{m-2r +1}_1(\Lambda V,{d})\not =0$ or  $H^{m}_1(\Lambda V,d)\not =0$. But, since $m>m-2r +1>m_1-2r +1=2r$ we conclude that $\dim H^{*}_1(\Lambda V,d)\geq 2$.
 This finishes the first step.
\subsubsection{Second step: $\dim H^*_2(\Lambda V,d)\geq 2$}
This step is the longest one and it will serve us to resume the general case, i.e. the third step.
In this step, we put $k=2$ in (\ref{Gysin}) to obtain  the following long exact sequence:
\begin{equation}\label{Gysin(2)} 
\ldots \stackrel{\delta ^*}{\rightarrow} H_{1}^{i-2r}(\Lambda V,d)\stackrel{j^*}{\rightarrow}H_{2}^{i}(\Lambda V,d)\stackrel{p^*}{\rightarrow}H_{2}^{i}(\Lambda W,\bar{d})\stackrel{\delta ^*}{\rightarrow}H_{2}^{i-2r+1}(\Lambda V,d)\rightarrow \ldots
\end{equation}
Noticing that $2r =|x_1| \leq   |x_2|\leq \ldots$, we obtain $H_2^{2r+1}(\Lambda V,d) =0$ so that, for $i=4r$, (\ref{Gysin(2)}) induces  the following exact sequence:
\begin{equation}\label{Gysin(2)(4r)}
0\rightarrow H_{1}^{4r-1}(\Lambda V,{d})\stackrel{p ^*}{\rightarrow}H_{1}^{4r-1}(\Lambda W,\bar{d})\stackrel{\delta ^*}{\rightarrow}H^{2r}_1(\Lambda V,d)\stackrel{j^*}{\rightarrow} H^{4r}_2(\Lambda V,d)\stackrel{p^*}{\rightarrow} H^{4r}_2(\Lambda W,\bar{d}){\rightarrow}0.
\end{equation}
It results that
\begin{equation}\label{(b)}
H^{4r}_2(\Lambda V,d)\cong \ker(p^*)\oplus H^{4r}_2(\Lambda W,\bar{d}).
\end{equation}
We need (as in the first step) to separate the case where $m_1 =2r$ from that where   $m_1 \not =2r$
\subsubsection*{Assume  $m_1=2r$} This implies that $H_1^{2r}(\Lambda V,d)\supseteq \mathbb{Q}.[x_1]\oplus \mathbb{Q}.[x_2]$. We should discuss three cases:
\\\\
{\bf (i)} If $\dim H^{4r}_2(\Lambda W,\bar{d})\geq 2$ then $\dim H^{*}_2(\Lambda V,{d})\geq 2$. 
\\\\
{\bf (ii)} If $\dim H^{4r}_2(\Lambda W,\bar{d})=1$ then, there are two sub-cases under consideration:
\\

{\bf (*)} In the first one, supposing $\ker (p^*)\not = 0$ (e.g.  $[x_1]^2\not = 0$ or $[x_1x_2]\not = 0$) we get $\dim H_2^{4r}(\Lambda V,d)\geq 2$ and consequently $\dim H_2^{*}(\Lambda V,d)\geq 2$.
\\

{\bf (**)} In the second one,  supposing  $\ker (p^*) = 0$, so,  using  (\ref{(b)}) 
we get the isomorphism 
$H_2^{4r}(\Lambda V,d)\cong H^{4r}_2(\Lambda W,\bar{d}) = \mathbb{Q}.[x_2]^2$.
 Therefore, by induction hypothesis,  we have $\dim H_2^*(\Lambda W,\bar{d})\geq 2$ and consequently \underline{there is some (least) integer $m>4r$} such that $H_2^{m}(\Lambda W,\bar{d})\not = 0$.
We continue by using the following exact sequence obtained from (\ref{Gysin(2)}) for  $i=m$:
\begin{equation}\label{Gysin(2)(m)}
\begin{array}{l}
\ldots \rightarrow H_{1}^{m-1}(\Lambda W,\bar{d})\stackrel{\delta ^*}{\rightarrow}H^{m-2r}_1(\Lambda V,d)\stackrel{j^*}{\rightarrow} H^{m}_2(\Lambda V,d)\stackrel{p^*}{\rightarrow}\\
\hspace*{5cm} H^{m}_2(\Lambda W,\bar{d}) 
\stackrel{\delta ^*}{\rightarrow}H^{m-2r+1}_2(\Lambda V,d){\rightarrow} \ldots 
\end{array}
\end{equation}
Clearly, If  $H^{m}_2(\Lambda V,d)\not =0$ then $\dim H^*_2(\Lambda V,d)\geq 2$.

 \underline{Next, we assume that $H^{m}_2(\Lambda V,d)=0$}. It results  that  $\delta ^*$ becomes a monomorphism and consequently $H^{m-2r+1}_2(\Lambda V,d)$ is non-zero so that $\dim H^*_2(\Lambda V,d)\geq 2$ {\bf unless}  if $m-2r+1=4r$ or equivalently, $m=6r-1$. \\
  \underline{We then proceed by assuming that $m=6r-1$} in which case  $\delta ^*$ is an isomorphism given by $\delta ^*([x_2x'_3])=[x_2]^2$ for some {unique} $x'_3\in V^{4r-1}$ such that $d(x'_3)=x_1x_2$. Indeed, if there exists another $x^{''}_3\in V^{4r-1}$ such that $\delta ^*([x_2x^{''}_3])=[x_2]^2$  and $d(x^{''}_3)=x_1x_2\not = 0$, we will have $d(x'_3-x^{''}_3)=0$, that is, there is in $V$ a cocycle with odd degree. This contradicts the hypothesis (a)  of the theorem. Similarly, there exists an unique $x'_4\in V^{4r-1}$ such that $d(x'_4)=x_1^2$. It results that $x_1x'_3-x_2x'_4$ is a  cocycle in $(\Lambda  V)^{6r-1=m}$.
 But,  since we assumed $H^{m}_2(\Lambda V,d)=0$ we must have $[x_1x'_3-x_2x'_4]=0$.  So, there exists $x'_5\in V^{6r-2}$ such that $d(x'_5)=x_1x'_3-x_2x'_4$ and consequently, $W'=W\backslash \{x_2\} $ satisfies $3\leq \dim W'\leq n-2$. In particular, $\dim H_2^*(\Lambda W', \bar{\bar{d}})\geq 2$ by induction hypothesis.
 
    We continue by
   considering the exact sequence:
    \begin{equation}\label{Gysin(2)(m')}
         \begin{array}{l}
         \ldots \rightarrow H_{1}^{m'-1}(\Lambda W',\bar{\bar{d}})\stackrel{\delta ^*}{\rightarrow}H^{m'-2r}_1(\Lambda W,\bar{d})\stackrel{j^*}{\rightarrow} H^{m'}_2(\Lambda W,\bar{d})\stackrel{p^*}{\rightarrow} \\
         \hspace*{6cm}H^{m'}_2(\Lambda W',\bar{\bar{d}})\stackrel{\delta ^*}{\rightarrow}H^{m}_2(\Lambda W,\bar{d}){\rightarrow} \ldots
               \end{array}
         \end{equation}
          obtained from (\ref{Gysin(2)}) with $(\Lambda W,\bar{d})$ (resp. $(\Lambda W',\bar{\bar{d}})$) replacing $(\Lambda V,d)$ (resp. $(\Lambda W,\bar{d})$) and  $i=m'=8r-2$, i.e. such that $m=m'-2r+1=6r-1$.
   Two  situations are  under consideration:\\
    $ \diamond$  If $H^{m'}_2(\Lambda W,\bar{d})\not = 0$, then, by reconsidering the exact sequence (\ref{Gysin(2)(m)}) with $m'$ instead of $m$ and noticing that $H^{m'-2r+1=m}_2(\Lambda V,d)=0$ (as it is assumed), we deduce that $H^{m'}_2(\Lambda V,{d})\not = 0$ and consequently $\dim H^*(\Lambda V,d)\geq 2$.
\\
   $ \diamond \diamond$ If $H^{m'}_2(\Lambda W,\bar{d})= 0$,  we have two possibilities:

   $\bullet$ { Firstly}, if $H^{m'}_2(\Lambda W',\bar{\bar{d}})\not = 0$
     then, the (last) morphism $\delta ^*$ in (\ref{Gysin(2)(m')}) being injective implies that $H^{m}_2(\Lambda W,\bar{d})\not = 0$, hence, using the isomorphism $\delta ^*$ in (\ref{Gysin(2)(m)}) we see that $\dim H^{m}_2(\Lambda W,{\bar{d}})=\dim H^{4r}(\Lambda V,d)=1$. Thus, $\dim H^{m'}_2(\Lambda W',\bar{\bar{d}})=1$. This assures the existence of an integer \underline{$m" \not = m'$} such that $H^{m"}_2(\Lambda W',\bar{\bar{d}})\not = 0$ by which once again (\ref{Gysin(2)(m')}) with $m"$ instead of $m'$ implies that either $H^{m"}_2(\Lambda W,\bar{d})\not =0$ or $H^{m"-2r+1}_2(\Lambda W,\bar{d})\not =0$. In particular, since $m$ is the least integer satisfying $H^{m}_2(\Lambda W,\bar{d})\not =0$, we deduce that $m">m=6r-1$ or, even more, $m"> m' = m +2r-1=8r-2$. Next, using (\ref{Gysin(2)(m)}) with $m">6r-1$ (resp. with  $m">8r-2$) instead of $m$, we get either $H^{m"}_2(\Lambda V,{d})\not =0$ or $H^{m"-2r+1}_2(\Lambda V,{d})\not =0$ (resp. $H^{m"-2r+1}_2(\Lambda V,{d})\not =0$ or $H^{m"-4r+2}_2(\Lambda V,{d})\not =0$). In all cases,  we obtain $\dim H^*(\Lambda V,d)\geq 2$.

 $\bullet \bullet$ { Secondly}, if $H^{m'}_2(\Lambda W',\bar{\bar{d}})= 0$, then,  once again by induction, 
  there is an integer $m">m'=8r-2$ such that $H^{m"}_2(\Lambda W',\bar{\bar{d}})\not = 0$. Thus, using one more  (\ref{Gysin(2)(m')}) with $m"$ instead of $m'$ we get $H^{m"}_2(\Lambda W,\bar{d})\not =0$ or $H^{m"-2r+1}_2(\Lambda W,\bar{d})\not =0$ and we conclude as in ($\bullet$) juste above.
\\\\
   {\bf (iii)} If $H^{4r}_2(\Lambda W,\bar{d})=0$, then the morphism $j^*: H_{1}^{2r}(\Lambda V,d)\stackrel{j^*}{\rightarrow}H_{2}^{4r}(\Lambda V,d))$ in (\ref{Gysin(2)(4r)}) is onto.
  Recall that (in general) the morphism $j^*: H_{k-1}^{i-2r}(\Lambda V,d)\stackrel{j^*}{\rightarrow}H_{k}^{i}(\Lambda V,d)$ in  (\ref{Gysin}) is defined as follows:
\begin{equation}\label{j}
\left\{ \begin{array}{l}
j^*([\chi'])=[x_1\chi']\\
\delta^*[\chi] = [\chi']
\end{array}
\right.
\;\;
    \hbox{with}\; \chi'\in \Lambda V, \; \chi \in \Lambda W \; \hbox{such that} \; d(\chi) = x_1\chi'.
\end{equation}
     In particular,  $[x_2^2]$ (as a class in $H_{2}^{4r}(\Lambda V,d)$) is necessarily zero, but  $[x_1]^2$ and $[x_1x_2]$
     may be non-zero. If they are both nonzero, then $\dim H_2^{*}(\Lambda V,d)\geq 2$. If $[x_1]^2\not =0$ and $[x_1x_2]=0$, we consider the equations $d(x'_3)=x_1x_2$ and  $d(x'_4)=x_2^2$  (some $x'_3$ and $x'_4$) which induce a new class $[x_2x'_3-x_1x'_4]\in H^{6r-1}_2(\Lambda V,d)$. If this class is non-zero we conclude that $\dim H_2^*(\Lambda V,d)\geq 2$. If not, it is still possible that  $H^{6r-1}_2(\Lambda V,d)\not =0$ which allows us to conclude immediately. 
     
     Now, if $H^{6r-1}_2(\Lambda V,d) =0$, and if moreover, there is another non-zero class, say $[x_3]\in H_1^{2r}(\Lambda V,d)$ (or even more than one) then,  $[x_2x_3] = [x_3]^2=0$ since they are not in $Im(j^*)$. But, $[x_1x_3]$ can be non-zero. If it is the case,  we conclude that $\dim H_2^{*}(\Lambda V,d)\geq 2$, otherwise  we may assume that 
     $H^{4r}_2(\Lambda V,d)\cong \mathbb{Q}.[x_1^2].$
     We then consider (\ref{Gysin(2)(m)}) with $m=6r-1$ by which   the morphism $\delta^*$ is actually an isomorphism and consequently, $\dim H^{6r-1}_2(\Lambda W,d)=1$. In particular, $m=6r-1$ is  effectively the least integer satisfying $H^{m}_2(\Lambda W,d)\not =0$. Using the induction hypothesis, we  introduce an $m'>m$ satisfying $H^{m'}_2(\Lambda W,d)\not =0$. Hence, making use once again of (\ref{Gysin(2)(m)}) with $m'$ instead of $m$ we have either  $H^{m'}_2(\Lambda V,d)\not = 0$ or else $H^{m'-2r+1}_2(\Lambda V,d)\not = 0$.  Since $m'-2r+1 > m-2r+1 = 4r$, we conclude, in both cases, that $\dim H^{*}_2(\Lambda V,d)\geq 1$.
     
    A similar argument may be used if $[x_1]^2=0$ and $[x_1x_2]\not =0$ by using instead the equations $d(x'_3)=x_1^2$ and  $d(x'_4)=x_2^2$, but this time, there is no new class to consider. So, we make use of the induction hypothesis to introduce an $m>4r$ such that $H^{m}_2(\Lambda W,d)\not =0$. We then conclude by using (\ref{Gysin(2)(m)}) which implies that $\dim H_2^m(\Lambda V,d) \oplus \dim H_2^{m-2r+1}(\Lambda V,d)\geq 1$. 
    
     It remains to discuss the last case, namely, when $[x_1]^2=[x_2]^2=[x_1x_2]=0$. So, based on the discussion just above, 
     we are in the following situation:
     $$H_2^{4r}(\Lambda V,d) = H_2^{4r}(\Lambda W,\bar{d})=0.$$
     To continue, we make use of the induction hypothesis on $H_2^{*}(\Lambda W,\bar{d})$ which give us an integer $m>4r$ (resp. two successive integers  $m'>m>4r$)   such that $\dim H_2^{m}(\Lambda W,\bar{d})\geq 2$  (resp. $\dim H_2^{m}(\Lambda W,\bar{d}) = 1$ and $\dim H_2^{m'}(\Lambda W,\bar{d})\geq 1$).
We have once again two possibilities:\\
$(*)$ Firstly, we assume  $\dim H_2^{m}(\Lambda W,\bar{d})\geq 2$.
     By the  exact sequence (\ref{Gysin(2)(m)}), we have necessarily $\dim H^{m}_2(\Lambda V,d)+\dim H^{m-2r+1}_2(\Lambda V,d)\geq 2$.\\ 
$(**)$ Secondly, we assume there are $m'>m>4r$ such that {$\dim H_2^{m}(\Lambda W,\bar{d}) = 1$} and {$\dim H_2^{m'}(\Lambda W,\bar{d})\geq 1$}. Notice that we may have $m'-2r+1 =m$.
Two sub-cases appear:
\\
       $\diamond$ If $H^{m-2r+1}_2(\Lambda V,d)\not = 0$, thus, since $m-2r+1<m<m'$ and  $H_2^{m'}(\Lambda W,\bar{d})\not = 0$, by using  (\ref{Gysin(2)(m)}) with $m'$ instead of $m$, we  obtain $H^{m'-2r+1}_2(\Lambda V,d)\not =0$ or else $H^{m'}_2(\Lambda V,d)\not =0$. Hence, in both cases, $\dim H^{*}_2(\Lambda V,d)\geq 2$.
\\
    $\diamond \diamond$  If $H^{m-2r+1}_2(\Lambda V,d)= 0$, from (\ref{Gysin(2)(m)}), we deduce that $H^{m}_2(\Lambda V,d)\not = 0$, so it remains to discuss the case where 
    $$\dim H^{m}_2(\Lambda V,d)=1.$$
    
   If $m\not = m'-2r+1$, we finish  just as  above  by re-using  (\ref{Gysin(2)(m)}) with $m'$ instead of $m$ and the fact that $\dim H_2^{m'}(\Lambda W,\bar{d})\geq 1$.
   
 Otherwise, \underline{$m= m'-2r+1$} so  $m$ and $m'$ have opposite parities. We are lead to consider the following relevant sub-cases:

  $\bullet$  In the first sub-case \underline{we assume $m$ even, hence $m'-2r$ is odd} and consequently, by hypothesis (a) of the theorem, we have $H^{m'-2r}_1(\Lambda V,d)=0$. Therefore,   once again, (\ref{Gysin(2)(m)}) with $m'$ instead of $m$,  give us the following exact sequence:
\begin{equation}\label{Gysin(m')(m even)}
    0\rightarrow H^{m'}_2(\Lambda V,d)\stackrel{p^*}{\rightarrow} H^{m'}_2(\Lambda W,\bar{d})\stackrel{\delta ^*}{\rightarrow}H^{m}_2(\Lambda V,d)\rightarrow \ldots 
\end{equation}
 Now, since $\dim H_2^{m}(\Lambda V,{d})=1$, we have either  $\dim H_2^{m'}(\Lambda W,\bar{d})\geq 2$  by which and  the exact sequence (\ref{Gysin(m')(m even)}) we  obtain $\dim H_2^{*}(\Lambda V,{d})\geq 2$, or else, {$\dim H_2^{m'}(\Lambda W,\bar{d})=1$} which implies that
     $\delta ^*$ is an isomorphism and consequently \underline{{$H^{m'}_2(\Lambda V,d)=0$}}.     
To continue, 
recall that \underline{until now, each of $H_2^{m}(\Lambda V,{d})$, $H_2^{m}(\Lambda W,\bar{d})$ and $H_2^{m'}(\Lambda W,\bar{d})$} \underline{ is one dimensional}.

Let then $[\xi] = [x_ix_j+x_{i'}x_{j'}+\ldots]\in H_2^{m}(\Lambda W,\bar{d})$  and $[\xi '] = [x_kx_l+x_{k'}x_{l'}+ \ldots]\in H_2^{m'}(\Lambda W,\bar{d})$,  some $j\geq i\geq 2$ and $l\geq k\geq 2$, be the respective basis elements. 
Since $m>4r$ we should have   $\mid x_i\mid>2r$ or $\mid x_j\mid>2r$, that is, if $i=2$ then $j>2$.  Similarly, since $m'>m>4r$ and $m'$ is odd, necessarily $\mid x_k\mid>2r$ or $\mid x_l\mid>2r$  and $l>k$. Clearly the same discussion holds  for $x_{i'}x_{j'}$ and $x_{k'}x_{l'}$. Thus, if there is more than one monôme in $[\xi]$ or in $[\xi ']$, we conclude that for $W' =W\backslash \{x_2\}$ we have $\dim W' \geq 3$. 
Otherwise, we have $[\xi] = [x_ix_j]$ and   $[\xi '] = [x_kx_l]$. In such a case,  if $i>2$ and $k>2$,  even if $k=i$ or $k=j$, we obtain $l\not = j$ due to the parities of $x_j$ and $x_l$. Once again,  $\dim W' \geq 3$. Now, if $k=2$, using the isomorphism 
$\delta^*$  we may put $\delta^*([\xi']) = [\alpha \xi]$ (some $\alpha \in \mathbb{Q}^*$), which implies $d(\xi') = \alpha \xi$. This gives $d(x_2x_l) = \alpha x_1x_ix_j$ which give us $|x_l|=|x_ix_j|=m$ is even, but this  contradics the fact that $m'=|x_kx_l|=|x_2x_l|$ is odd. Next, we assume  that $i=2$. Using once again the isomorphism 
$\delta^*$, we get $d(x_kx_l) = \alpha x_1x_2x_j=\alpha(x_1x_2)x_j= \alpha d(x_tx_j)$ (where $d(x_t) = x_1x_2$) or, equivalently, $d(x_kx_l-\alpha x_jx_t)=0$. But,  since $|x_kx_l-\alpha x_jx_t|=m'$ and $H_2^{m'}(\Lambda V,{d})=0$, there exist $x_s\in V$ such that $d(x_s) = x_kx_l-\alpha x_jx_t$. In particular, even if $k=j$ (resp. $l=j$), we have $\{x_j, x_l, x_s\}\subseteq W'$ (resp. $\{x_k, x_l, x_s\}\subseteq W'$) and consequently, $3\leq \dim W'\leq n-2$. 
         It results, from the discussion above, that we may use the sequence (\ref{Gysin(2)(m')}) with the acctual subspace $W'$. Therefore, a similar discussion to that made in the sub-case $(ii)-(**)$  above, {especially after assuming that $m=6r-1$},  give us the required conclusion by using (\ref{Gysin(2)(m)})  with $m'$ (which is odd) instead of $m=6r-1$  and (\ref{Gysin(2)(m')}) with $m"=m'+2r-1$ (which is even) instead of $m'=8r-2$.

  $\bullet \bullet$ In the second sub-case,  \underline{ we  assume that  $m$ is odd and $m'$ is even}. Recall  that $m=m'-2r+1>4r$, $\dim H_2^{m}(\Lambda W,\bar{d}) = 1$,  $\dim H_2^{m'}(\Lambda W,\bar{d})\geq 1$ and  $\dim H^{m}_2(\Lambda V,d)=1$.
  
  $\triangleright$ If $\dim H_2^{m'}(\Lambda W,\bar{d})\geq 2$, to conclude,  it suffice to use the sequence (\ref{Gysin(2)(m)}) with $m'$ instead of $m$ by which we obtain $\dim H^{m'}_2(\Lambda V,d)\geq 1$. 
  
  $\triangleright \triangleright$ If $\dim H_2^{m'}(\Lambda W,\bar{d}) =1$, in  (\ref{Gysin(2)(m)}) with $m'$ instead of $m$ we still have $\delta^*$ is an isomorphism (of one dimensional spaces), so, $j ^*: H_1 ^{m'-2r}(\Lambda V,d)\rightarrow  H_2^{m'}(\Lambda V,{d})$ is an epimorphism. Moreover, $m'-2r$ being  actually even, by hypothesis, we may have $H_1 ^{m'-2r}(\Lambda V,d)\not =0$. If this is the case, $H_2^{m'}(\Lambda V,{d}) \not =0$  and we are done. If not i.e. if $H_2^{m'}(\Lambda V,{d})=0$, then,
  after swaping rules of $m$ and $m'$, the discussion made in $(\bullet)$ just above remains truth and we are also done.

\subsubsection{Suppose $m_1>2r$:}
  Therefore, \underline{$H_{2}^{4r}(\Lambda W,\bar{d})=0$ and $\dim H^{4r}_2(\Lambda V,d)=0$}. Indeed,  since the least order of a (possible) non-zero class in $H_{2}^{*}(\Lambda W,\bar{d})$ should be  $2m_1$, then $H_{2}^{4r}(\Lambda W,\bar{d})=0$. Next, since $H_1^{2r}(\Lambda V,d) =\mathbb{Q}.[x_1]$, using
    (\ref{Gysin(2)(4r)})  we obtain $H^{4r}_2(\Lambda V,d)= \ker(p^*) = Im(j^*)$ is generated by   $j^*[x_1]=[x_1^2]$. Now, (\ref{j}) implies that   
$d(\chi)=x_1^2$  
for some $\chi \in \Lambda W^{4r-1}$.
Therefore, {$\dim H^{4r}_2(\Lambda V,d)=0$}. 

By considering (\ref{Gysin(2)}) with $i=|x_1| +|x_2|=m_1+2r$, we get the exact sequence:
 \begin{equation}\label{(Gysin(m_1+2r))}
  \ldots {\rightarrow} H^{m_1+2r-1}_1(\Lambda V,{d})\stackrel{p^*}{\rightarrow} H^{m_1+2r-1}_1(\Lambda W,\bar{d})\stackrel{\delta ^*}{\rightarrow} H^{m_1}_1(\Lambda V,d)\stackrel{j^*}{\rightarrow} H^{m_1 +2r}_2(\Lambda V,d){\rightarrow} 0.
  \end{equation}
  Notice that  $m_1+2r< 2|x_2|$ implies $H^{m_1+2r}_2(\Lambda W,\bar{d})=0$. Here, as in the previous discussion, we have three cases:
  \\\\
  {\bf (i)}  if $\dim  H^{m_1 +2r}_2(\Lambda V,d)\geq 2$, we finish.
  \\\\
  {\bf (ii)} If $\dim  H^{m_1 +2r}_2(\Lambda V,d) =1$ then, by (\ref{(Gysin(m_1+2r))}), $H^{m_1}_1(\Lambda V,d)\not =0$. Therefore, by hypothesis (a) of the theorem,
   \underline{ $m_1$ should be even} and  $H^{m_1+2r-1}_1(\Lambda V,{d})=0$. The exact sequence (\ref{(Gysin(m_1+2r))}) induces then the following one:
 \begin{equation}\label{Gysin(m_1(qcq))}
 0 {\rightarrow} H^{m_1+2r-1}_1(\Lambda W,\bar{d})\stackrel{\delta ^*}{\rightarrow} H^{m_1}_1(\Lambda V,d)\stackrel{j^*}{\rightarrow} H^{m_1 +2r}_2(\Lambda V,d){\rightarrow} 0.
 \end{equation}
We continue by considering once again  (\ref{Gysin}) with $i=2m_1$. This gives the exact sequence:
 \begin{equation}\label{(Gysin(2m_1))}
 \ldots {\rightarrow}H^{2m_1-1}_1(\Lambda W,\bar{d})\stackrel{\delta ^*}{\rightarrow}  H^{2m_1-2r}_1(\Lambda V,d)\stackrel{j^*}{\rightarrow} H^{2m_1}_2(\Lambda V,d)\stackrel{p^*}{\rightarrow} H^{2m_1}_2(\Lambda W,\bar{d}){\rightarrow}0.
 \end{equation}
In fact, $H^{2m_1-2r+1}_2(\Lambda V,d)=0$ since  $2m_1-2r+1<2m_1$ and $2m_1-2r+1 \not = m_1+2r$ (here $m_1$ is even). There are two sub-cases under consideration:

 {\bf (*)} In the first one, assuming $H^{2m_1}_2(\Lambda W,\bar{d})\not =0$  we get $H^{2m_1}_2(\Lambda V,d)\not =0$ and we are done.

 {\bf (**)} In the second one, assuming $H^{2m_1}_2(\Lambda W,\bar{d}) =0$, it is still possible that $H^{2m_1}_2(\Lambda V,d)\not =0$. If that is the case, we are done. If no, we use  the induction hypothesis for $(\Lambda W, \bar{d})$ to get an integer  $m>2m_1$ (which we assume the smallest one) such that \underline{$H^{m}_2(\Lambda W,\bar{d})\not =0$}. Remark that, since $[x_1^2]=[x_2^2]=0$, there must exist $x'_3,\; x'_4\in W$ such that $d(x'_3)=x_1^2$ and $d(x'_4)=x_2^2$ so that effectively $\dim W\geq 3$ as required by the inductive hypothesis.
 Two situations are under consideration: 

 $\diamond$ Suppose firstly that $m_1>4r-1$, so that $m-2r+1>m_1+2r$. Hence, if $H_2^{m-2r+1}(\Lambda V,d)\not =0$ we have $\dim H_2^{*}(\Lambda V,d)\geq 2$. But, \underline{if $H_2^{m-2r+1}(\Lambda V,d) =0$} then, by considering once again  (\ref{Gysin}) with $i=m$ we  get a copy of (\ref{(Gysin(2m_1))}) {with $m$ instead of  $2m_1$} which implies  that  $H^{m}_2(\Lambda V,d)\stackrel{p^*}{\rightarrow} H^{m}_2(\Lambda W,\bar{d})$ is onto. It results that   $H^{m}_2(\Lambda V,d)\not = 0$ and then $\dim H_2^{*}(\Lambda V,d)\geq 2$.

 $\diamond \diamond$ Suppose secondly that $m_1<4r-1$ and notice that    $m-2r+1=m_1+2r$ if and only if  $m=m_1+4r-1$. Thus,  as  $2m_1<m_1+4r-1$, we may have effectively $m-2r+1=m_1+2r$. Therefore,
 if $\dim H^{m}_2(\Lambda W,\bar{d})\geq 2$, using  again  the exact sequence (\ref{Gysin(2)(m)}),  we have necessarily  $\dim H_2^{*}(\Lambda V,d)\geq 2$.
But, if \underline{$\dim H^{m}_2(\Lambda W,\bar{d})= 1$}, we have either $m-2r+1\not =m_1+2r$ in which case $\dim H_2^{m}(\Lambda V,d) + \dim H_2^{m-2r+1}(\Lambda V,d)\geq 1$ and consequently $\dim H^*_2(\Lambda V,d)\geq 2$;  or else, $m-2r+1  = m_1+2r$ in which case we 
once again make use of  induction hypothesis to introduce an $m'>m$ for which $\dim H^{m'}_2(\Lambda W,\bar{d})\not = 1$. We conclude after  re-using (\ref{Gysin(2)(m)}) with $m'$ instead of $m$ and noticing that actually $m'-2r+1\not =m_1+2r$.
 \\\\
{\bf (iii)} Assume that \underline{$H^{m_1+2r}_2(\Lambda V,d)=0$}. So, since $H^{4r}_2(\Lambda V,d)=0$,  $2m_1-2r+1 = m_1 + (m_1-2r+1)<2m_1$,  and $m_1$ and $m_1-2r+1$ have distinct parities, we deduce, using hypothesis $(a)$ of the theorem, that $H^{2m_1-2r+1}_2(\Lambda V,d)=0$. Thus, if we put $2m_1$ instead of $m$ in (\ref{Gysin(2)(m)})  we obtain again the exact sequence (\ref{(Gysin(2m_1))}):
 \begin{equation*}\label{Gysin(2m_1)(2)}
 \ldots {\rightarrow}H^{2m_1-1}_1(\Lambda W,\bar{d})\stackrel{\delta ^*}{\rightarrow}   H^{2m_1-2r}_1(\Lambda V,d)\stackrel{j^*}{\rightarrow} H^{2m_1}_2(\Lambda V,d)\stackrel{p^*}{\rightarrow} H^{2m_1}_2(\Lambda W,\bar{d})\rightarrow 0.
  \end{equation*}
Clearly, if $\dim H^{2m_1}_2(\Lambda W,\bar{d})\geq 2$ we finish. It remains to discuss two other sub-cases:
  
{\bf (*)}  We first assume that $\dim H_2^{2m_1}(\Lambda W,\bar{d})=1$. Then, by the above exact sequence, we have $\dim H_2^{2m_1}(\Lambda V,d)\geq 1$. If this is greater than two, we are done, but, in case where \underline{$\dim H_2^{2m_1}(\Lambda V,d)=1$}, we use induction hypothesis to introduce some $m>2m_1$ such that $H_2^{m}(\Lambda W,\bar{d})\not =0$. We then make use of the  sequence (\ref{Gysin(2)(m)}) which we repeat here for convenience
\begin{equation*}
\begin{array}{l}
\ldots \rightarrow H_{1}^{m-1}(\Lambda W,\bar{d})\stackrel{\delta ^*}{\rightarrow}H^{m-2r}_1(\Lambda V,d)\stackrel{j^*}{\rightarrow} H^{m}_2(\Lambda V,d)\stackrel{p^*}{\rightarrow}\\
\hspace*{5cm} H^{m}_2(\Lambda W,\bar{d}) 
\stackrel{\delta ^*}{\rightarrow}H^{m-2r+1}_2(\Lambda V,d){\rightarrow} \ldots.
\end{array}
\end{equation*}
If   $\dim H^{m}_2(\Lambda W,\bar{d}) \geq 2$, then    $\dim H^{m}_2(\Lambda V,d) + \dim H^{m-2r+1}_2(\Lambda V,d)\geq 2$ and we are done. If  \underline{$\dim H^{m}_2(\Lambda W,\bar{d}) =1$}, we are also done if moreover $H^{m}_2(\Lambda V,d)\not =0$. But,  if  \underline{$H^{m}_2(\Lambda V,d)=0$} we deduce that  $\dim H^{m-2r+1}_2(\Lambda V,d)\geq 1$ and we are also done,  unless if \underline{$\dim H^{m-2r+1}_2(\Lambda V,d) =1$} and {$m-2r+1 =2m_1$}.
 At this stage, we have $m=2m_1+2r-1$  odd and
$$H^{m}_2(\Lambda V,d)=0,\; \dim H^{2m_1}_2(\Lambda V,d) = \dim H_2^{2m_1}(\Lambda W,\bar{d})= \dim H^{m}_2(\Lambda W,\bar{d}) =1.$$
We then finish, by introducing a subspace $W'$ such that  $W=\mathbb{Q}x_2\oplus W'$ using a discussion similar to that made is the sub-case $(m_1=2r)-(iii)-(\bullet)$ when $m$ (which correspond here to $2m_1$) is even and $m'$ (which correspond here to $2m_1+2r-1$) is odd.

 (**) Now, if $H^{2m_1}_2(\Lambda W,\bar{d})=0$, then $ H^{2m_1-2r}_1(\Lambda V,d)\stackrel{j^*}{\rightarrow} H^{2m_1}_2(\Lambda V,d)$ in the exact sequence (\ref{(Gysin(2m_1))}) is onto. 
 
 Since $2m_1-2r$ is even, by hypothesis $(a)$ of the Theorem \ref{Thm1}, we have three cases:  \\
 $\triangleright$ If $\dim H^{2m_1}_2(\Lambda V,d)\geq 2$ we clearly finish.\\
 $\triangleright \triangleright$ If $\dim H^{2m_1}_2(\Lambda V,d)=1$, then, using the inductive hypothesis,  we introduce an integer $m>2m_1$ such that $H_2^{m}(\Lambda W,\bar{d})\not = 0$ by which and (\ref{Gysin(2)(m)}) we also conclude unless when $H_2^{m}(\Lambda V,{d})=0$, $\dim H_2^{m}(\Lambda W,\bar{d}) =1$ and $m-2r+1 =2m_1$. We procede, using hypothesis on $(\Lambda W,\bar{d})$  to introduce an $m'>m$ such that $\dim H_2^{m'}(\Lambda W,\bar{d}) \not =0$. By means of (\ref{Gysin(2)(m)}) with $m'$ instead of $m$, the case requiring discution is when $\dim H_2^{m'}(\Lambda W,\bar{d}) =1$. But in such a case, we finich by noticing that either $\dim H_2^{m'}(\Lambda V,{d}) =1$ or else $\dim H_2^{m'-2r+1}(\Lambda V,{d}) =1$ with $m'-2r+1>2m_1$.\\
 $\triangleright \triangleright \triangleright$ If $\dim H^{2m_1}_2(\Lambda V,d)=0$. In this case, we use the inductive hypothesis to introduce $m'>m>2m_1$ such that $\dim H_2^{m}(\Lambda W,\bar{d})=1$ and $H_2^{m'}(\Lambda W,\bar{d})\geq 1$ and make use of (\ref{Gysin(2)(m)}). The extreme case is when $H_2^{m-2r+1}(\Lambda V,{d}) =0$, each of  $H_2^{m}(\Lambda V,{d})$,  $H_2^{m}(\Lambda W,\bar{d})$ and $H_2^{m'}(\Lambda W,\bar{d})$ is one dimensional, and $m =m'-2r+1$. We then proceed as in sub-case  $(m_1=2r)-(iii)-(\bullet)$, or else, sub-case  $(m_1=2r)-(iii)-(\bullet \bullet)$) according to the parity of $m$.

 
   This completes the proof for $k=2$.
 \subsubsection{General step: $\dim H_k^*(\Lambda V,d)\geq 2$}
  We assume that $\dim H_{k-1}^*(\Lambda V,d)\geq 2$ for $k-1\geq 2$.
 Notice first that by putting in the exact sequence (\ref{Gysin})   $i=2kr$ we obtain the following one:
 \begin{equation}\label{Gysing(k)}
 \ldots \rightarrow H_{k-1}^{2kr-1}(\Lambda W,\bar{d})\stackrel{\delta ^*}{\rightarrow} H_{k-1}^{2(k-1)r}(\Lambda V,d)\stackrel{j^*}{\rightarrow}H^{2kr}_{k}(\Lambda V,d)\stackrel{p^*}{\rightarrow}H^{2kr}_{k}(\Lambda W,\bar{d}){\rightarrow} 0.
 \end{equation}
 For, using the $1$-connectedness of $X$ we have $r\geq 2$, then, noticing that the degree of $k$ classes is at least $2kr$, we obtain $H_{k}^{2(k-1)r+1}(\Lambda V,d) =0$. It follows that
 \begin{equation}\label{(k)}
 H_k^{2kr}(\Lambda V,d) \cong \ker (p^*)\oplus H_k^{2kr}(\Lambda W,\bar{d}).
 \end{equation}
 {\it In this general case, we will repose our induction on $\dim H_k^{2kr}(\Lambda W,\bar{d})$ instead of $m_1$.
  Three cases depending on $\dim H_k^{2kr}(\Lambda W,\bar{d})$ are required.}
 \\\\
 {\it First,  assume that $\dim H_k^{2kr}(\Lambda W,\bar{d})\geq 2$.} 
  Thus, (\ref{(k)}) implies  that $\dim H_k^{*}(\Lambda V,d)\geq 2$.
 \\\\
 {\it Second, assume that $\dim H_k^{2kr}(\Lambda W,\bar{d})=1$.}   We necessarily have $m_1=2r$ and two sub-cases are under consideration:
 
 (*) In the first, we suppose  $\dim \ker (p^*)\geq 1$, so we also obtain $\dim H_k^{*}(\Lambda V,d)\geq 2$.
 
 (**) In the second, we suppose $\ker (p^*)=0$, then, by  (\ref{(k)}), we have $H_k^{2kr}(\Lambda V,d)\cong H_k^{2kr}(\Lambda W,\bar{d})\cong \mathbb{Q}.[x_2]^{k}$.
  By using  a similar discussion as that made in $\S 4.2.2 (ii) (**)$, we obtain  an integer $m> 2kr$ such that $H_{k}^m(\Lambda W, \bar{d})\not =0$ and a consequent exact sequence (which is similar to (\ref{Gysin(2)(m)})):
 \begin{equation}\label{Gysin(k)(m)}
 \begin{array}{l}
 \ldots \rightarrow H_{k-1}^{m-1}(\Lambda W,\bar{d})\stackrel{\delta ^*}{\rightarrow}H^{m-2r}_{k-1}(\Lambda V,d)\stackrel{j^*}{\rightarrow}  H^{m}_k(\Lambda V,d)\stackrel{p^*}{\rightarrow} H^{m}_k(\Lambda W,\bar{d}) \stackrel{\delta ^*}{\rightarrow} \\
  \hspace*{6cm} H^{m-2r+1}_k(\Lambda V,d){\rightarrow}\ldots 
 \end{array}
 \end{equation}
 Once again, we have
 
  If $H^{m}_k(\Lambda V,d)\not = 0$ then, $\dim H^{*}_k(\Lambda V,d)\geq 2$.
 
  If $H^{m}_k(\Lambda V,d)=0$,  the morphism $H^{m}_k(\Lambda W,\bar{d}) \stackrel{\delta ^*}{\rightarrow} H^{m-2r+1}_k(\Lambda V,d)$ becomes a monomorphism. Thus, $H^{m-2r+1}_k(\Lambda V,d)\not =0$ and $\dim H^{*}_k(\Lambda V,d)\geq 2$ except if $m-2r+1=2kr$ or equivalently $m=2(k+1)r-1$. 
  
  We then procced by assuming $H^{m}_k(\Lambda V,d)=0$ and $m=2(k+1)r-1$.
  It results that  $\delta ^*$ is an isomorphism given by $\delta ^*([(x_2)^{k-1}x'_3])=[x_2]^{k}$ for some $x'_3\in V^{4r-1}$ such that $d(x'_3)=x_1x_2$. Furthermore, we have also $[x_1]^{k}=0$ (since $\ker(p^*)=0$).
   Thus, there exists $x'_4$
    such that $d(x'_4)=x_1^2$
    so that $x_1^k=d(x_1^{k-2}x'_4)$. Therefore, using the same reasoning as in $\S 4.2.2 (ii) (**)$, we show that $x'_3$ is unique and there is some $x'_5$ saltsfying $d(x'_5) = x_1x'_3-x_2x'_4$. 
     We then get $W' = W\backslash \{x_2\}$ satisfying $3\leq \dim W' \leq n-2$. We next make use of
    the following exact sequence
   \begin{equation}\label{Gysin(k)(m')}
   \begin{array}{l}
   \ldots \rightarrow H_{k-1}^{m'-1}(\Lambda W',\bar{\bar{d}})\stackrel{\delta ^*}{\rightarrow}H^{m'-2r}_{k-1}(\Lambda W,\bar{d})\stackrel{j^*}{\rightarrow} H^{m'}_k(\Lambda W,\bar{d})\stackrel{p^*}{\rightarrow} \\
   \hspace*{6cm}  H^{m'}_k(\Lambda W',\bar{\bar{d}})\stackrel{\delta ^*}{\rightarrow} H^{m}_k(\Lambda W,\bar{d}){\rightarrow} \ldots
   \end{array}
  \end{equation}
  obtained with $m'=m+2r-1$. Notice that this later is similar to the sequence (\ref{Gysin(2)(m')}) already considered in the second step.
  To conclude this case, it suffices to use the same discussion as the one made in  the cases $\diamond$ and $\diamond \diamond$ of $\S 4.2.2 (ii) (**)$.
 \\\\
 {\it Third, assume that $H_k^{2kr}(\Lambda W,\bar{d})=0$.}  Thus, $H^{2(k-1)r}_{k-1}(\Lambda V,{d})\stackrel{j^*}{\rightarrow} H^{2kr}_k(\Lambda V,{d})$ in (\ref{Gysing(k)}) is onto. We should discuss two relevant sub-cases:
 
 $(*)$ {\it Assume $m_1 = 2r$}. We have to use a discussion similar to that made just after the equation (\ref{j}), but this time replacing respectively $[x_1^2]$, $[x_2^2]$ and $[x_1x_2]$ by $[x_1^k]$, $[x_2^k]$ and $[x_1^sx_2^t]$ (some paire or eventually paires of integers $s>0,\; t>0$ such that $s+t=k$). Here, we have $[x_2^k]=0$ as a class in $H_k^{2kr}(\Lambda V,{d})$ (since it is not in the image of $j^*$), but $[x_1^k]$ and the $[x_1^sx_2^t]'s$ may be non-zero. If at least two between such classes are non-zero, we are done. Otherwise, we proceed as follows:
 
 First, assume that  $[x_1^k]=0$ and only one of the $[x_1^sx_2^t]$'s is non-zero. Thus,  there exists $x'_4\in V^{6r-1}$ with $dx'_4=x_1^2$ so that $x_1^k = d(x_1^{k-2}x'_4)$. We put similarly, $dx'_3=x_2^2$, some $x'_3\in V^{6r-1}$, so that $x_2^k = d(x_2^{k-2}x'_3)$. The relations $dx'_3=x_2^2$ and $dx'_4=x_1^2$ do not induce an eventual cocycle so, as in $\S 4.2.2 (iii)$ we make use of induction hypothesis to introduce some $m>2kr$ such that $H^m_2(\Lambda W,\bar{d})\not = 0$. We then conclude by using (\ref{Gysin(k)(m)}) which ensures   that $\dim H^m_2(\Lambda V,\bar{d}) \oplus H^{m-2r+1}_2(\Lambda V,{d})\geq 1$.
 
 Second, assume that all the $x_1^sx_2^t$'s are coboudaries and $[x_1^k]\not =0$. Then, 
  there exists, say an $x'_4\in V^{6r-1}$ such that $d(x'_4)=x_1x_2$ which give us $x_1^sx_2^t = d(x'_4x_1^{s-1}x_2^{t-1})$. Using moreover  $dx'_3=x_2^2$  we obtain a cocycle $x_1x'_3-x_2x'_4$ wich may define a non-zero class $[x_1^{s-1}x_2^{t-1}(x_1x'_3 - x_2x'_4)]\in H^{2(k+1)r-1}(\Lambda V,d)$. If it is effectively non-zero, we finish, if not, we use again the inductive hypothesis to introduce some $m>2kr$ and conclude just as above.
 
 Notice here that, by definition of $j^*$ and since it is onto, if there is anothe non-zero class of the form $[x_3^{k-1}]\in H_1^{2(k-1)r}(\Lambda V,d)$ (or even more than one) then, we should have $[x_3^k]=0$.
 
 It remain then to discuss \underline{the situation where $H_k^{2kr}(\Lambda V,d)= H_k^{2kr}(\Lambda W,\bar{d})=0$}. 
 Specifically, we use induction hypothesis on $H_k^{*}(\Lambda W,\bar{d})$ to introduce either (a): some $m>2kr$ such that $\dim H_k^{m}(\Lambda W,\bar{d})\geq 2$ or else (b): two integers $m>2kr$ and $m'>m>2kr$, with the possibility that $m=m'-2r+1$, such that $\dim H_k^{m}(\Lambda W,\bar{d})=1$ and $\dim H_k^{m'}(\Lambda W,\bar{d})\geq 1$.
 
 The sub-case requiring special treatment is (b)  when \underline{$m=m'-2r+1$} and
  \underline{$\dim H_k^{m}(\Lambda W,\bar{d})=1$, $H_k^{m'}(\Lambda W,\bar{d})\not = 0$,  $H_k^{m-2r+1}(\Lambda V,d)=0$}. This  situation is similar to $\S 4.2.2-(iii)-(**)(\diamond \diamond)$ where $m=m'-2r+1$. Clearly, $m$ and $m'$ have opposite parities.
  We should then consider the following relevant sub-cases:
  
  $\diamond$ Assume $m$ even and $m'$ odd. By $(\ref{Gysin(k)(m)})$, we have $\dim H_k^m(\Lambda V,d)\not =0$. Hence, if $\dim H_k^m(\Lambda V,d)\geq 2$, we finish, if not, we have $\dim H_k^m(\Lambda V,d)=1$ so, $p^*:  H^{m}_k(\Lambda V,d)\stackrel{p^*}{\rightarrow} H^{m}_k(\Lambda W,\bar{d})$ becomes an isomorphism, in particulr, $Im(j^*)=0$. Therefore, we have either $\dim H_k^{m'}(\Lambda W,\bar{d})\geq 2$ in which case, by $(\ref{Gysin(k)(m)})$ with $m'$ instead of $m$, we deduce that $\dim H_k^{m'}(\Lambda V,d)\geq 1$ and we are done; otherwise, \underline{$\dim H_k^{m'}(\Lambda W,\bar{d})=1$}. In particular the morphism $\delta^*$ in $(\ref{Gysin(k)(m)})$ becomes an isomorphism (and it is the case when we take $'$ instead of $m$).
  
  We then continue by introdicing the generating classes $[\xi] = [x_{i_1}x_{i_2}x_{i_3}\ldots x_{i_k} + \ldots]$ and $[\xi'] = [x_{j_1}x_{j_2}x_{j_3}\ldots x_{j_k} + \ldots]$  of $H_k^{m}(\Lambda W,\bar{d})$ and $H_k^{m'}(\Lambda W,\bar{d})$ respectively, where $i_k\geq \ldots \geq i_2\geq i_1\geq 2$ and $j_k\geq \ldots \geq j_2\geq j_1\geq 2$. In a similar way to the case where $k=2$, we introduce $W' = W\backslash \{x_2\}$ satisfying $3\leq \dim W' \leq n-2$ and translate to this level, the discussion made in the sub-case $(**)$ of $(ii)$  \underline{especially  when $m=6r-1$}. This   give us the required conclusion by using (\ref{Gysin(k)(m)})  with $m'$ (which is odd) instead of $m=6r-1$  and (\ref{Gysin(k)(m')}) with $m"=m'+2r-1$ (which is even) instead of $m'=8r-2$.
  
 $\diamond \diamond$ Now, we assume $m$ odd and $m'$ even. The discussion made at the beginning of the sub-case $\diamond$ just above remains available and we have to treat the case where $\dim H_k^{m}(\Lambda W,\bar{d})=\dim H_k^{m'}(\Lambda W,\bar{d})= \dim H_k^{m}(\Lambda V,d)=1$. Once again, sweeping the roles of $m$ and $m'$ we obtain an adequate $W' = W\backslash \{x_2\}$  satisfying $3\leq \dim W' \leq n-2$ by which we reach the same conclusion as just above. This finishes the case where $m_1=2r$.
  
  (**) {\it Assume $m_1>2r$}. In particular, we have $H^{2kr}_{k}(\Lambda W,\bar{d})=0$, thus, using  (\ref{Gysing(k)}), we deduce that $H_k^{2kr}(\Lambda V,d)\cong Im(j^*)=\ker(p^*)$. We have two cases:
  
  $\bullet$ {\it If $\dim H_k^{2kr}(\Lambda V,d) =1$} (by hypothesis, $H^{2(k-1)r}_{k-1}(\Lambda W,\bar{d})$ may be non-zero).
  Thus,   inspired by the discussion made in the second step, we put in (\ref{Gysin}) $i=m_1+2(k-1)r$. It results in the exact sequence:
   \begin{equation}
     \ldots \stackrel{p^*}{\rightarrow} H^{m_1+2(k-1)r-1}_{k-1}(\Lambda W,\bar{d})\stackrel{\delta ^*}{\rightarrow} H^{m_1+2(k-2)r}_{k-1}(\Lambda V,d)\stackrel{j^*}{\rightarrow} H^{m_1 +2(k-1)r}_k(\Lambda V,d){\rightarrow} 0.
     \end{equation}
     Here again, $H^{m_1 +2(k-1)r}_k(\Lambda W,\bar{d})=0$, since the  least degree of a cocycle with length $k$ is $m_1 +2(k-1)r<km_1$. We must distinguish between two sub-cases
 
    $\diamond$ In the first, we assume  that $H^{m_1 +2(k-1)r}_k(\Lambda V,d)\not =0$, so, $\dim  H^{*}_k(\Lambda V,d)\geq 2$.
 
    $\diamond \diamond$  In the second, we assume that $H^{m_1 +2(k-1)r}_k(\Lambda V,d)=0$.
 We then  consider the following exact sequence obtained from (\ref{Gysin}) with $i=2m_1+2(k-2)r > 2kr$:
   \begin{equation}
    \ldots {\rightarrow}H^{2m_1+2(k-2)r-1}_{k-1}(\Lambda W,\bar{d})\stackrel{\delta ^*}{\rightarrow}  H^{2m_1+2(k-3)r}_{k-1}(\Lambda V,d) \stackrel{j^*}{\rightarrow} H^{2m_1+2(k-2)r}_k(\Lambda V,d){\rightarrow}0
    \end{equation}
  ($j^*$ is onto since $2m_1+2(k-2)r<km_1$). Thus,  we are led to use the inductive assumption:
  $$H^{jm_1 +2(k-j)r}_k(\Lambda V,d)=0;\; \; \forall\;  1\leq j\leq k-2.$$
  Once again, the exact sequence obtained from (\ref{Gysin}) with $i=(j+1)m_1+2(k-j-1)r$ permits to conclude (at this stage) that,
 
   $\circ$ either $H^{2m_1+2(k-2)r}_k(\Lambda V,d)\not =0$, hence,  $\dim  H^{*}_k(\Lambda V,d)\geq 2$ or else,
 
  $\circ \circ$ $H^{2m_1+2(k-2)r}_k(\Lambda V,d)=0$, hence,  $H^{jm_1 +2(k-j)r}_k(\Lambda V,d)=0$ for $1\leq j\leq k-1$. That is, we are in a situation where:
  $$\dim H^{2kr}_k(\Lambda V,d)=1,\; H^{jm_1 +2(k-j)r}_k(\Lambda V,d)=0; \; \forall 1\leq j\leq k-1.$$
 
  We then argue by considering the following ecaxt sequence obtained from (\ref{Gysin}) when $i=km_1$:
 \begin{equation}\label{Gysin(k)(km_1)}
 \begin{array}{l}
 \ldots {\rightarrow}H^{km_1-1}_{k-1}(\Lambda W,\bar{d})\stackrel{\delta ^*}{\rightarrow}  H^{km_1-2r}_{k-1}(\Lambda V,d)\stackrel{j^*}{\rightarrow} H^{km_1}_k(\Lambda V,d)\\
 \hspace*{6cm} \stackrel{p^*}{\rightarrow} H^{km_1}_k(\Lambda W,\bar{d})\stackrel{\delta ^*}{\rightarrow}H^{km_1-2r+1}_k(\Lambda V,d)\rightarrow \ldots
 \end{array}
  \end{equation}
  Notice neverthless that,  as $km_1-2r+1  > 2(k-1)r + 1$, we may have $km_1-2r+1 = 2kr$ since then $m_1 = 2r+\frac{2r-1}{k}$ which my be a positive integer.  Two possibilities are under consideration:
 
  ($\star$) if $H^{km_1}_k(\Lambda W,\bar{d})\not =0$, then   $H^{km_1}_k(\Lambda V,d)\not =0$ or   $H^{km_1-2r+1}_k(\Lambda V,d)\not =0$. In both cases   $\dim  H^{*}_k(\Lambda V,d)\geq 2$ unless when  $\dim H^{km_1}_k(\Lambda W,\bar{d})=1$, $H^{km_1}_k(\Lambda V,d)=0$, and $km_1-2r+1=2kr$. In such a case, we make use of the induction hypothesis to introduce an $m>km_1$ satisfying $H^{m}_k(\Lambda W,\bar{d})\not =0$. So, by reconsidering the exact sequence (\ref{Gysin(k)(m)})  and noticing  that,  necessarily, $m-2r+1> 2kr$, we conclude that $\dim  H^{*}_k(\Lambda V,d)\geq 2$.
 
  ($\star \star$) If $H^{km_1}_k(\Lambda W,\bar{d})=0$, we proceede by  considering the following exact sequence obtained from (\ref{Gysin}) with $i=km_1+2r$:
     \begin{equation}
     \begin{array}{l}
      \ldots {\rightarrow}H^{km_1+2r-1}_{k-1}(\Lambda W,\bar{d})\stackrel{\delta ^*}{\rightarrow}  H^{km_1}_{k-1}(\Lambda V,d)\stackrel{j^*}{\rightarrow} H^{km_1+2r}_k(\Lambda V,d)\rightarrow \\
      \hspace*{6cm} \stackrel{p^*}{\rightarrow} H^{km_1+2r}_k(\Lambda W,\bar{d})\stackrel{\delta ^*}{\rightarrow}H^{km_1+1}_{k}(\Lambda V,d)\rightarrow \ldots
      \end{array}
      \end{equation} 
      Thus, if $H^{km_1+2r}_k(\Lambda W,\bar{d})\not = 0$ then $\dim H^{km_1+2r}_k(\Lambda V,d) +  \dim H^{km_1+1}_k(\Lambda V,d)\geq 1$ and we are done. If not, 
   we use the induction hypothesis to introduce some $m> km_1+2r$ such that $H^{m}_k(\Lambda W,\bar{d})\not =0$. Hence, usince once again (\ref{Gysin(k)(m)})     we have either $H^{m}_k(\Lambda V,d)\not = 0$ or else $H^{m-2r+1}_k(\Lambda V,d)\not =0$, therefore, since $m-2r+1 > km_1+1 > 2kr$, we are also done.
 
  $\bullet \bullet$  {\it If $H_k^{2kr}(\Lambda V,d)=0$}. Since  $H_k^{2kr}(\Lambda W,\bar{d})=0$, we are in a situation similar to that of $\S 4.2.3$.
  
   The induction hypothesis on $H_k^{*}(\Lambda W,\bar{d})$ enables us to consider a least integer $m>2kr$ (resp. two least integers  $m'>m>2kr$),   such that $\dim H_k^{m}(\Lambda W,\bar{d})\geq 2$  (resp. $\dim H_k^{m}(\Lambda W,\bar{d}) = 1$ and $H_k^{m'}(\Lambda W,\bar{d})\not =0$). That is we have to discuss the following sub-cases:
 
   $\diamond$ In the first, we assume  $\dim H_k^{m}(\Lambda W,\bar{d})\geq 2$. It results from the exact sequence (\ref{Gysin(k)(m)})  that   $$\dim H_k^{m-2r+1}(\Lambda V,d) +\dim H_k^{m}(\Lambda V,d)\geq 2$$ and thus $\dim H_k^{*}(\Lambda V,d)\geq 2$. 
        
   $\diamond \diamond$ In the second, we assume that $\dim H_k^{m}(\Lambda W,\bar{d}) = 1$ and $H_k^{m'}(\Lambda W,\bar{d})\not =0$ so, by using  (\ref{Gysin(k)(m)}), we have either (i) $H_k^{m-2r+1}(\Lambda V,d)\not =0$ or else (ii) $H_k^{m}(\Lambda V,d)\not =0$. In the case (i), by using  once more (\ref{Gysin(k)(m)}) with $m'$ instead of $m$, we get either $H^{m'}_k(\Lambda V,d)\not =0$ or  else $H_k^{m'-2r+1}(\Lambda V,d)\not =0$. In both cases, we finish since    $m'>m'-2r+1> m-2r+1$. In the  case (ii), that is when $H_k^{m-2r+1}(\Lambda V,d) =0$, the case requiring discussion is when, $H^{m'}(\Lambda V,d) =0$, $\dim H^m_k(\Lambda V,d)= \dim H^{m'-2r+1}_k(\Lambda V,d) = \dim H_k^{m'}(\Lambda W,\bar{d} = 1$. Here, we conclude provided that $m\not =m'-2r+1$.  Now, if $m=m'-2r+1$, we finish by using a similar discussion as in the sub-cases $(m_1=2r)-(*)-(\diamond)$ and $(m_1=2r)-(*)-(\diamond \diamond) $  just above.

 This completes this general case and consequently the proof of Theorem 1.3.

\end{document}